\documentclass{article}
\usepackage{spconf}
\usepackage{graphicx}
\usepackage{amsmath}
\usepackage{amsfonts}
\usepackage{amssymb}
\usepackage{epsfig}
\usepackage{algorithm}
\usepackage{algorithmic}
\usepackage{url}

\newcommand{\bo}[1]{{\bf #1}}

\title{Total Variation Restoration of Speckled Images\\ Using a Split-Bregman Algorithm}
\name{Jos\'{e} M. Bioucas-Dias \hspace{1.5cm} M\'{a}rio A. T. Figueiredo}
\address{{\it Instituto de Telecomunica\c{c}\~{o}es}, \\
{\it Instituto Superior T\'{e}cnico}, \\
Lisboa, {\bf Portugal}\\
 Email: $\{$jose.bioucas,\ mario.figueiredo$\}$@lx.it.pt}

\ninept

\begin{document}

\maketitle

\begin{abstract}
Multiplicative noise models occur in the study of several
coherent imaging systems, such as synthetic aperture radar and sonar, and
ultrasound and laser imaging. This type of noise is also commonly referred
to as {\it speckle}. Multiplicative noise introduces two additional
layers of difficulties with respect to the popular
Gaussian additive noise model: (1) the noise is multiplied by (rather than added to)
the original image, and (2) the noise is not Gaussian, with Rayleigh and Gamma being commonly
used densities. These two features of the multiplicative noise model preclude the
direct application of state-of-the-art restoration methods, such as those based
on the combination of total variation or wavelet-based regularization with a
quadratic observation term. In this paper, we tackle these difficulties by:
(1) using the common trick of converting the multiplicative model into
an additive one by taking logarithms, and (2) adopting the recently proposed
split Bregman approach to estimate the underlying image under total variation
regularization. This approach is based on formulating a constrained problem
equivalent to the original unconstrained one, which is then solved
using Bregman iterations (equivalently, an augmented Lagrangian method).
A set of experiments show that the proposed method yields state-of-the-art results.
\end{abstract}

%the power
%read  by the sensor from a given
%resolution element  is  the contributions of many coherent signals. Given
%that these signals differ from a random phase, they interfere randomly in a destructive or constructive manner, leading to what is known as speckle noise or multiplicative noise.
%These observations follow, very often, a Rayleigh distribution
%(or Gamma, if multilook/averaging is used);


\section{Introduction}
\subsection{Coherent Imaging and Speckle Noise}
The standard statistical models of coherent imaging systems, such as synthetic aperture radar/sonar (SAR/SAS),
ultrasound imaging,  and laser imaging, are supported on multiplicative noise mechanisms.
With respect to a given resolution cell of the imaging device, a coherent system acquires the so-called
in-phase and quadrature components\footnote{The in-phase and quadrature components are the
outputs of two demodulators with respect to, respectively,  $\cos(\omega_0 t)$) and $\sin(\omega_0 t)$,
where $\omega_0$ is the carrier angular frequency.} which are collected in a complex reflectivity (with
the in-phase and quadrature components corresponding  to the real and imaginary parts, respectively).
The complex reflectivity of a given resolution cell results from the contributions of
all the individual scatterers present in that cell, which interfere in a destructive or constructive manner,
according to their spatial configuration. When this configuration is random, it yields random
fluctuations of the complex reflectivity, a phenomenon which is termed {\em speckle}.
The statistical properties of speckle have been widely studied and there is a large body of literature  \cite{conf:Goodman:JOSA:76}, \cite{Oliver}.
Assuming no strong specular reflectors and a large number of randomly distributed  scatterers in each resolution cell (relative to the carrier wavelength),
the squared amplitude ({\em intensity}) of the complex reflectivity
is exponentially distributed \cite{Oliver}.
The term {\it multiplicative noise} is clear from the following
observation: an exponential  random variable can be written as the product of its mean
value (parameter of interest) by an exponential variable of unit mean (noise).  The scenario
just described,  known as {\em fully developed speckle}, leads to observed intensity images with a
characteristic granular appearance due to the very low {\em signal to noise ratio} (SNR).
Notice that the SNR, defined as the ratio between the squared intensity mean and the
intensity variance, is equal to one ($0\,$dB).

\subsection{Restoration of Speckled Images: Previous Work}
A common approach to improving the SNR in coherent imaging consists in
averaging independent observations of the same pixel. In SAR/SAS systems,
this procedure is called {\em multi-look} ($M$-look, in the case of $M$ looks),
and each independent observation may be  obtained by a different segment of the sensor array.
For fully developed speckle, the SNR of an  $M$-look
image is $M$. Another way to obtain an M-look image is to low pass filter (with a moving
average kernel with support size $M$) a 1-look fully developed speckle image, making  evident the
tradeoff between SNR and spatial resolution. A great deal of research has been devoted to
developing nonuniform filters which average large numbers of pixels in homogeneous regions
yet avoid smoothing across discontinuities in order to preserve image detail/edges \cite{art:FrostStilesPAMI82}.
Many other speckle reduction techniques have been proposed; see \cite{Oliver} for a comprehensive
literature review up to 1998.

A common assumption is that the underlying reflectivity image is
piecewise smooth. In image restoration under multiplicative noise, this assumption has been
formalized using  Markov random fields, under the Bayesian framework \cite{Oliver}, \cite{Bioucas-Dias98} and, more recently, using {\em total variation} (TV) regularization  \cite{AubertAujol08}, \cite{Huang09}, \cite{RudinLionsOsher03}, \cite{ShiOsher07}.

\subsection{Contribution}
In this paper, we adopt TV regularization. In comparison with the canonical
additive Gaussian noise  model, we face  two difficulties: the noise is multiplicative;
the noise is non-Gaussian, but follows Rayleigh or Gamma distributions.
We tackle these difficulties by first converting the multiplicative model
into an  additive  one (which is a common procedure) and then adopting the
recently proposed split Bregman approach to solve the optimization problem
that results from adopting a total variation regularization criterion.

Other works that have very recently addressed the restoration of speckled images
using TV regularization include \cite{AubertAujol08}, \cite{Huang09},
\cite{RudinLionsOsher03}, \cite{ShiOsher07}. The commonalities and differences
between our approach and the ones followed in those papers will be
discussed after the detailed description of our method, since this discussion
requires notation and concepts which will be introduced in the next section.

\section{Problem Formulation}
Let $\bo{y}\in \mathbb{R}_{+}^{n}$ denote an $n$-pixels observed image,
assumed to be a sample of a random image $\bo{ Y}$, the mean of which is the
underlying reflectivity image $\bo{x}\in \mathbb{R}_{+}^{n}$, {\it i.e.}, $\mathbb{E}[\bo{Y}] = \bo{x}$.
Adopting a conditionally independent multiplicative noise model, we have
\begin{equation}
Y_i = x_i N_i ,\;\; \mbox{for\  $\; i=1,...,n$,}\label{eq:multiply}
\end{equation}
where ${\bf N}\in \mathbb{R}_{+}^{n}$ is an image of independent and identically distributed (iid)
noise random variables with unit mean, $\mathbb{E}(N_i)=1$, following a common density $p_N$.
For $M$-look fully developed speckle noise, $p_N$ is a Gamma density with $E[N]=1$, and $\sigma_N^2=1/M$,
{\it i.e.},
\begin{equation}
   \label{eq:Gamma}
   p_N(n) = \frac{M^M}{\Gamma(M)}\; n^{M-1}e^{-nM}.
\end{equation}

An additive noise model is obtained by taking logarithms of (\ref{eq:multiply}).
For some pixel of the image, the observation model becomes
\begin{equation}
      \label{eq:log_observation}
       \underbrace{\log Y}_G = \underbrace{\log x}_z + \underbrace{\log N}_W.
\end{equation}
The density of the random variable $W = \log N$ is
\begin{equation}
      p_W(w) = p_N(n=e^w)\,e^w =  \frac{M^M}{\Gamma(M)}\; e^{Mw}e^{-e^w M},
\end{equation}
thus
\begin{equation}
    p_{G|Z}(g|z) = p_W(g-z).
\end{equation}

Under the regularization and Bayesian frameworks, the original image is inferred by
solving a minimization problem with the form
\begin{equation}
       \widehat{\bo{z}}   \in \arg\min_{\bo{z}} L(\bo{z}),  \label{eq:L_uncons}
\end{equation}
where $L(\bo{z})$ is the penalized minus log-likelihood,
\begin{eqnarray}
L(\bo{z}) & = & -\log p_{\bo{G}|\bo{Z}}(\bo{g}|\bo{z}) +\lambda\, \phi(\bo{z}).\\
& = & M \sum_{s=1}^n  \left(z_s + e^{g_s-z_s}\right) + \lambda\, \phi(\bo{z}) + A,\label{eq:neg_like}
\end{eqnarray}
with $A$ an irrelevant additive constant,  $\phi$ the
penalty/regularizer (negative of the log-prior, from  a the Bayesian perspective),
and $\lambda$ the regularization parameter.

In this  work, we adopt the TV regularizer, that is,
\begin{equation}
 \phi(\bo{z})  =  \mbox{TV}(\bo{z}) = \sum_{s=1}^n \sqrt{(\Delta^h_s\bo{z})^2+(\Delta^v_s\bo{z})^2},\label{eq:theTV}
\end{equation}
where  $(\Delta^h_s\bo{z}$ and $\Delta^v_s\bo{z})$ denote the horizontal and vertical first
order differences at pixel $s\in\{1,\dots,n\}$, respectively.

Each term $\left(z_s + e^{g_s-z_s} \right)$ of (\ref{eq:neg_like}), corresponding to the negative log-likelihood, is strictly convex and coercive, thus so is their sum.
Since the TV regularizer is also convex (though not strictly so),
the objective function $L$ possesses a unique minimizer \cite{CombettesSIAM}.  In terms
of optimization, these are desirable properties that would not hold if we had
formulated the inference in the original variables $\bo{x}$, since the resulting
negative log-likelihood is not convex; this was the approach followed
in \cite{AubertAujol08} and \cite{RudinLionsOsher03}.

\section{Bregman/Augmented Lagrangian  Approach}

There are several efficient algorithms to compute the TV regularized
solution for a quadratic data term; for recent work, see \cite{Chambolle04},
\cite{ChanGolubMulet}, \cite{conf:Mario:BJN:ICIP:00}, \cite{GoldsteinOsher},
\cite{WangYangYinZhang}, \cite{ZhuWrightChan}, and other references therein.
When the data term is not quadratic, as in (\ref{eq:neg_like}),
the problem is more difficult and far less studied. Herein, we follow the split
Bregman approach \cite{GoldsteinOsher} which is composed of the following two
steps: (splitting) a constrained problem equivalent to the original unconstrained one
is formulated; (Bregman) this constrained problem is solved using the Bregman iterative approach \cite{YinOsherGoldfarbDarbon}. Before describing these two steps in detail, we briefly review the Bregman iterative approach.
The reader is referred to \cite{GoldsteinOsher}, \cite{YinOsherGoldfarbDarbon}, for more details.

\subsection{Bregman Iterations}
Consider a constrained optimization problem of the form
\begin{eqnarray}
\nonumber \min_{\bf x} & & E({\bf x})\\
\mbox{s.t.} & &H({\bf x}) = 0,\label{eq:constrained_general}
\end{eqnarray}
with $E$ and $H$ convex, $H$ differentiable, and $\min_{\bf x}H({\bf x}) = 0$ . The so-called  Bregman divergence associated with the convex function $E$ is defined as
\begin{equation}
    D_{E}^{\mathbf{p}}(\mathbf{x},\mathbf{y}) \equiv E(\mathbf{x})-E(\mathbf{y}) - \langle\mathbf{p},\mathbf{x}-\mathbf{y}\rangle,
\end{equation}
where $\mathbf{p}$ belongs to the subgradient of $E$ at ${\bf y}$, {\it i.e.},
%({\it i.e.,} $\mathbf{p}\in\partial E(\mathbf{y})$).
$$\mathbf{p}\ \in\ \partial E(\mathbf{y})=\{\mathbf{u}:E(\mathbf{x})\geq E(\mathbf{y})+\langle\mathbf{u},\mathbf{x}-\mathbf{y}\rangle, \ \forall \mathbf{x} \in \mathbf{dom} E \}. $$

The Bregman iteration is given by
\begin{eqnarray}
 {\bf x}^{k+1} &=& \arg\min_{\bf x}  D_{E}^{\mathbf{p}^k}(\mathbf{x},\mathbf{x}^{k}) +  H({\bf x}) \label{eq:updatex} \\
&=& \arg\min_{\bf x} E({\bf x}) - \langle\mathbf{p}^k,\mathbf{x}-\mathbf{x}^k\rangle +  H({\bf x}),\label{eq:updatex2}
\end{eqnarray}
where $\mathbf{p}^k\in\partial E(\mathbf{x}^k)$. It has been shown that this procedure
converges to a solution of (\ref{eq:constrained_general}) \cite{GoldsteinOsher},\cite{YinOsherGoldfarbDarbon}.

Concerning the update of $\mathbf{p}^k$,  we have from \eqref{eq:updatex}, that $\mathbf{0}\in\partial (D_{E}^{\mathbf{p}}(\mathbf{x},\mathbf{x}^{k}) + \mu H({\bf x}))$, when this sub-differential is evaluated at $\mathbf{x}^{k+1}$, that is
$$\mathbf{0}\in\partial (D_{E}^{\mathbf{p}}(\mathbf{x}^{k+1},\mathbf{x}^{k}) +  H({\bf x}^{k+1})).$$
Since it was assumed that $H$ is differentiable, and since $\mathbf{p}^{k+1}\in \partial E(\mathbf{x}^{k})$ at this point, $\mathbf{p}^{k+1}$ should be chosen as
\begin{equation}
    \mathbf{p}^{k+1} = \mathbf{p}^{k} - \nabla H(\mathbf{x}^{k+1}).
\end{equation}

In the particular case where $H(\bo{x}) = (\tau/2)\|\bo{Ax-b}\|_2^2$, it can be shown (see \cite{YinOsherGoldfarbDarbon}) that
the iteration (\ref{eq:updatex2}) is equivalent to
\begin{eqnarray}
    {\bf x}^{k+1} &=& \arg\min_{\bf x} E(\bo{x}) + \frac{\tau}{2} \|\bo{Ax - b}^k\|_2^2 \label{eq:Bregman1a}\\
    \bo{b}^{k+1} & = & \bo{b} + \bo{b}^k - \bo{A}\bo{x}^k .\label{eq:Bregman1b}
\end{eqnarray}

\subsection{Splitting the Problem into a Constrained Formulation}
The original unconstrained problem (\ref{eq:L_uncons}) is equivalent to the constrained formulation
 \begin{eqnarray}
      \label{eq:L_cons}
      (\widehat{\bo{z}},\widehat{\bo{u}}) & = &
\arg\min_{\bo{z},\bo{u}} L(\bo{z},\bo{u})\\
      \label{eq:z_u_c}
      \mbox{s.t.} & &      \|\bo{z}-\bo{u}\|_2^2 = 0,
\end{eqnarray}
 with
\begin{eqnarray}
      \label{eq:L_cons_obj}
     L(\bo{z},\bo{u})  = M \sum_{s=1}^n \left( z_s + e^{g_s-z_s}\right) + \lambda\,\mbox{TV}(\bo{u}).
\end{eqnarray}

Notice how the original variable (image) $\bo{z}$ is split into a pair of variables $(\bo{z,u})$, which are decoupled in the objective function (\ref{eq:L_cons_obj}).

\subsection{Applying Bregman Iterations}
Notice that the problem  (\ref{eq:L_cons})-(\ref{eq:z_u_c}) has exactly the
form (\ref{eq:constrained_general}), with ${\bf x} \equiv [\bo{z}^T \bo{u}^T]^T$,
$E({\bf x}) \equiv L(\bo{z,u})$, and $H({\bf x}) \equiv (\tau/2) \|{\bf A x}-{\bf b}\|_2^2$, with ${\bf A} = [{\bf I}, -{\bf I}]$
and ${\bf b} = \bo{0}$. Using this equivalence, the Bregman iteration (\ref{eq:Bregman1a})-(\ref{eq:Bregman1b}) becomes
\begin{eqnarray}
     \label{eq:Breg_simpler_iter}
    (\bo{z}^{k+1}, \bo{u}^{k+1}) & = & \arg\min_{\bo{z},\bo{u}}
L(\bo{z},\bo{u})+\frac{\tau}{2}\| \bo{z}-\bo{u}-\bo{b}^k\|^2,\label{eq:Bregman_u_z}\\
    \bo{b}^{k+1}  & = &  \bo{b}^{k} - (\bo{z}^{k} -\bo{u}^{k}).
\end{eqnarray}

We address the minimization in (\ref{eq:Bregman_u_z}) using an alternating
minimization scheme with respect to ${\bf u}$ and ${\bf z}$. The complete resulting
algorithm is summarized in Algorithm 1.

\begin{algorithm}
\label{alg:} \caption{TV restoration of multilook images.}
\begin{algorithmic}[1]
\REQUIRE $\bo{z}=0$, $\bo{u}=0$, $\bo{b}=0$, $\lambda$, $\tau$, $k := 1$.
 \REPEAT
  \FOR{$t=1:t_m$}
   \STATE $\bo{z}^k := \arg\min_{\bo{z}} \sum_{s=1}^n \left( z_s +
e^{g_s-z_s}\right)+\frac{\tau}{2M}\| \bo{z}\!-\!\bo{u}^k\!-\!\bo{b}^k\|^2$
   \STATE $\bo{u}^k := \arg\min_{\bo{u}} \frac{1}{2}\|
\bo{u}-\bo{z}^k+\bo{b}^k\|^2+\frac{\lambda}{\tau}\,\mbox{TV}(\bo{u})$.
   \ENDFOR
   \STATE  $\bo{b}^{k+1}   :=   \bo{b}^{k} - (\bo{z}^{k} -\bo{u}^{k})$
   \STATE $k := k + 1$
   \UNTIL{$\|{\bf z}^{k}-{\bf z}^{k-1}\|_2^2/\|{\bf z}^{k-1}\|_2^2 < 10^{-4}$}
 \end{algorithmic}
 \end{algorithm}

The minimization with respect to $\bo{z}$, in line 3, has closed form
in terms of the Lambert W function \cite{Corless}. However, we found
that the Newton method yields a faster solution by running just four iterations.
Notice that the minimization in line 3 is in fact a set of $n$ decoupled
scalar minimizations. For the minimization with respect to $\bo{u}$ (line 4),
which is a TV denoising problem, we run  a few iterations (typically 10) of
Chambolle's algorithm \cite{Chambolle04}. The number of inner iterations $t_m$
was set to one in all the experiments reported below. The stopping criterion
(line 8) is the same as in \cite{Huang09}. The estimate of ${\bf x}$ produced
by the algorithm is naturally $\widehat{\bf x} = e^{{\bf z}^k}$, component-wise.

Notice how the split Bregman approach converted a difficult problem
involving a non-quadratic term and a TV regularizer into two simpler
problems: a decoupled minimization problem (line 3) and a TV denoising
problem with a quadratic data term (line 4).

\subsection{Remarks}
In the case of linear constraints, the Bregman iterative procedure
defined in (\ref{eq:updatex2}) is equivalent to an augmented Lagrangian
method \cite{Nocedal}; see \cite{TaiWu}, \cite{YinOsherGoldfarbDarbon} for proofs.
It is known that the augmented Lagrangian is better conditioned that the standard Lagrangian
for the same problem, thus a better numerical behavior is expectable.

TV-based image restoration under multiplicative noise was recently addressed in
\cite{ShiOsher07}. The authors apply an inverse scale space flow, which
converges to the solution of the constrained problem of minimizing $\mbox{TV}({\bf z})$
under an equality constraint on the log-likelihood; this requires a carefully
chosen stopping criterion, because the solution of this constrained problem
is not a good estimate.

In \cite{Huang09}, a splitting of the variable is also used to obtain
an objective function with the form
\begin{equation}
E({\bf z,u}) = L({\bf z,u}) + \alpha \|{\bf z - u}\|_2^2;
\end{equation}
this is the so-called splitting-and-penalty method. Notice that the
minimizers of $E({\bf z,u})$ converge to those of (\ref{eq:L_cons})-(\ref{eq:z_u_c})
only when $\alpha$ approaches infinity. However, since
$E({\bf z,u})$ becomes severely ill-conditioned when $\alpha$ is very
large, causing numerical difficulties, it is only practical to minimize
$E({\bf z,u})$ with moderate values of $\alpha$; consequently, the solutions obtained
are not minima of the regularized negative log-likelihood (\ref{eq:neg_like}).

\section{Experiments}
In this section we report experimental results comparing the
performance of the proposed approach
with that of the recent state-of-the-art methods in \cite{AubertAujol08}
and \cite{Huang09}.
All the experiments use synthetic data, in the sense that the observed image is
generated according to (\ref{eq:multiply})-(\ref{eq:Gamma}), where ${\bf x}$ is a
clean image. As in \cite{Huang09}, we select the regularization parameter
$\lambda$ by searching for the value leading to the lowest mean squared error
with respect to the true image. The algorithm is initialized with the
observed noisy image. The quality of the estimates is assessed using
the relative error  (as in \cite{Huang09}),
\[
\mbox{Err} = \frac{\|\widehat{\bf x} - {\bf x}\|_2}{\|{\bf x}\|_2}.
\]

Table~\ref{tab:results} reports the results obtained using Lena and
the Cameraman as original images, for the same values of the
number of looks ($M$ in (\ref{eq:Gamma})) as used in \cite{Huang09}.
In these experiments, our method always achieves lower relative errors
with fewer iterations, when compared with the  methods from \cite{Huang09}
and \cite{AubertAujol08} (the results concerning
the algorithm from \cite{AubertAujol08} are those reported in \cite{Huang09}).
It's important to point out that the computational cost of each iteration of the
algorithm of \cite{Huang09} is essentially the same as that of our algorithm.

\begin{table}
\centering
\caption{Experimental results (Iter denotes the number of iterations; Cam. is the
Cameraman image).}\label{tab:results}
\vspace{0.3cm}
\begin{tabular}{l l | l l | l l | l l }
\hline\hline
      &     & \multicolumn{2}{c|}{Proposed} & \multicolumn{2}{c|}{\cite{Huang09}} & \multicolumn{2}{c}{\cite{AubertAujol08}}\\
Image & $M$  & Err   & Iter & Err    & Iter & Err   & Iter \\ \hline
{\footnotesize Lena} & 5    &  {\footnotesize 0.1134} & 53   & {\footnotesize 0.1180} & 115  & {\footnotesize 0.1334} & 652   \\
{\footnotesize Lena} & 33   &  {\footnotesize 0.0688} & 23   & {\footnotesize 0.0709} & 178  & {\footnotesize 0.0748} & 379   \\
{\footnotesize Cam.} & 3    &  {\footnotesize 0.1331} & 100  & {\footnotesize 0.1507} & 182  & {\footnotesize 0.1875} & 1340   \\
{\footnotesize Cam.} & 13   &  {\footnotesize 0.0892} & 97   & {\footnotesize 0.0989} & 196 &  {\footnotesize 0.1079} & 950   \\
\hline
\end{tabular}
\end{table}

\begin{figure}
\centering
\includegraphics[width=4cm]{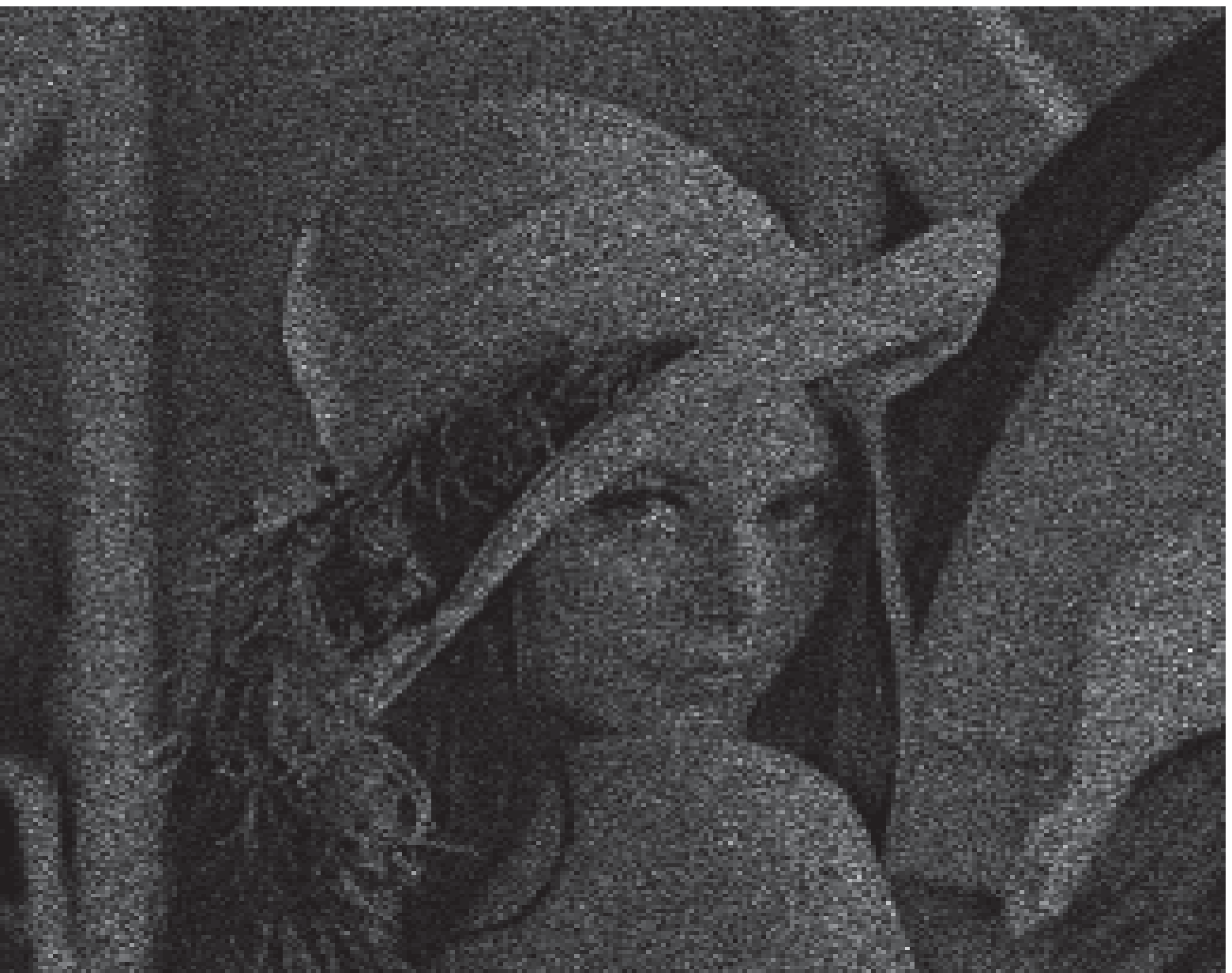}\hspace{0.1cm}\includegraphics[width=4cm]{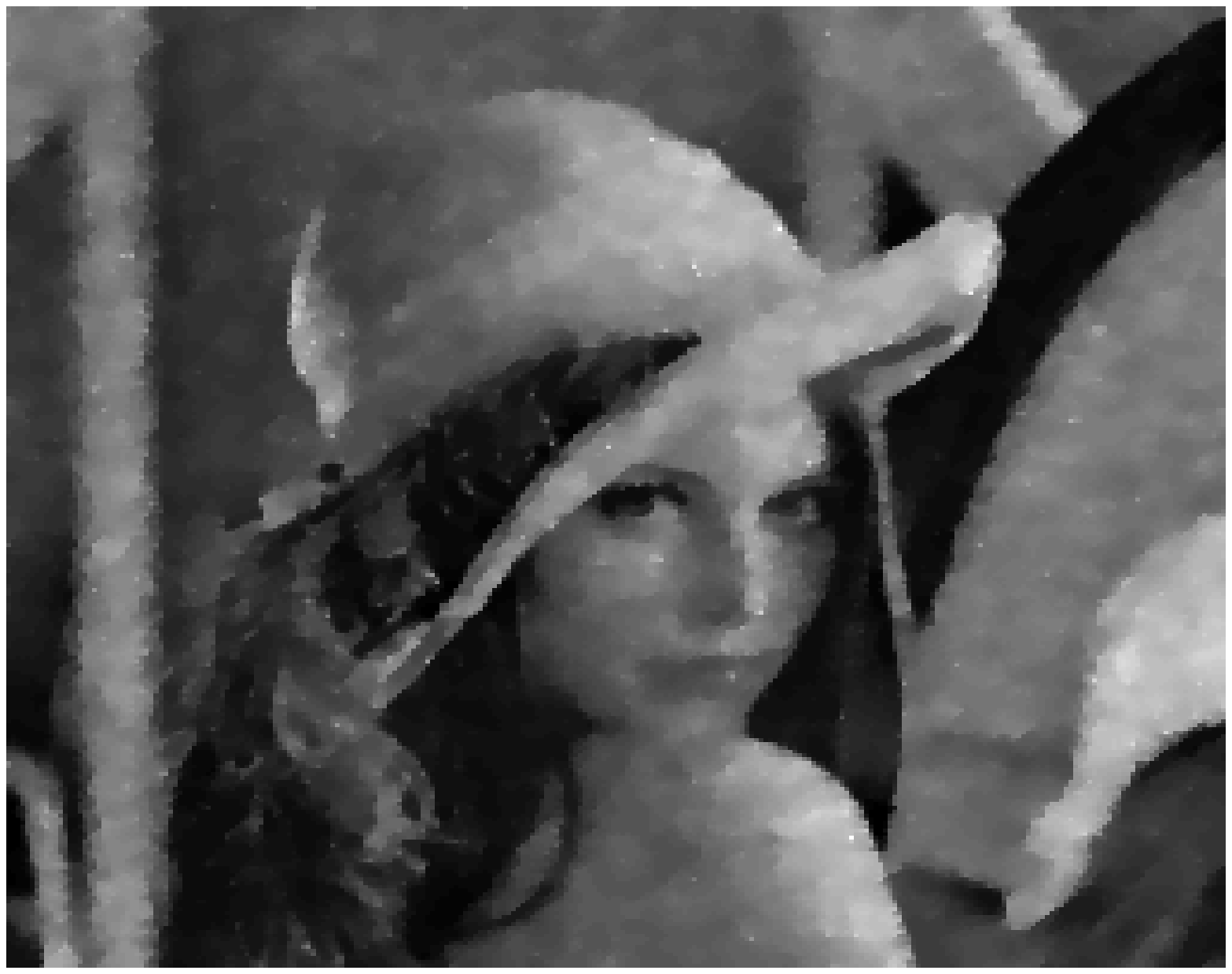}

\vspace{0.1cm}
\includegraphics[width=4cm]{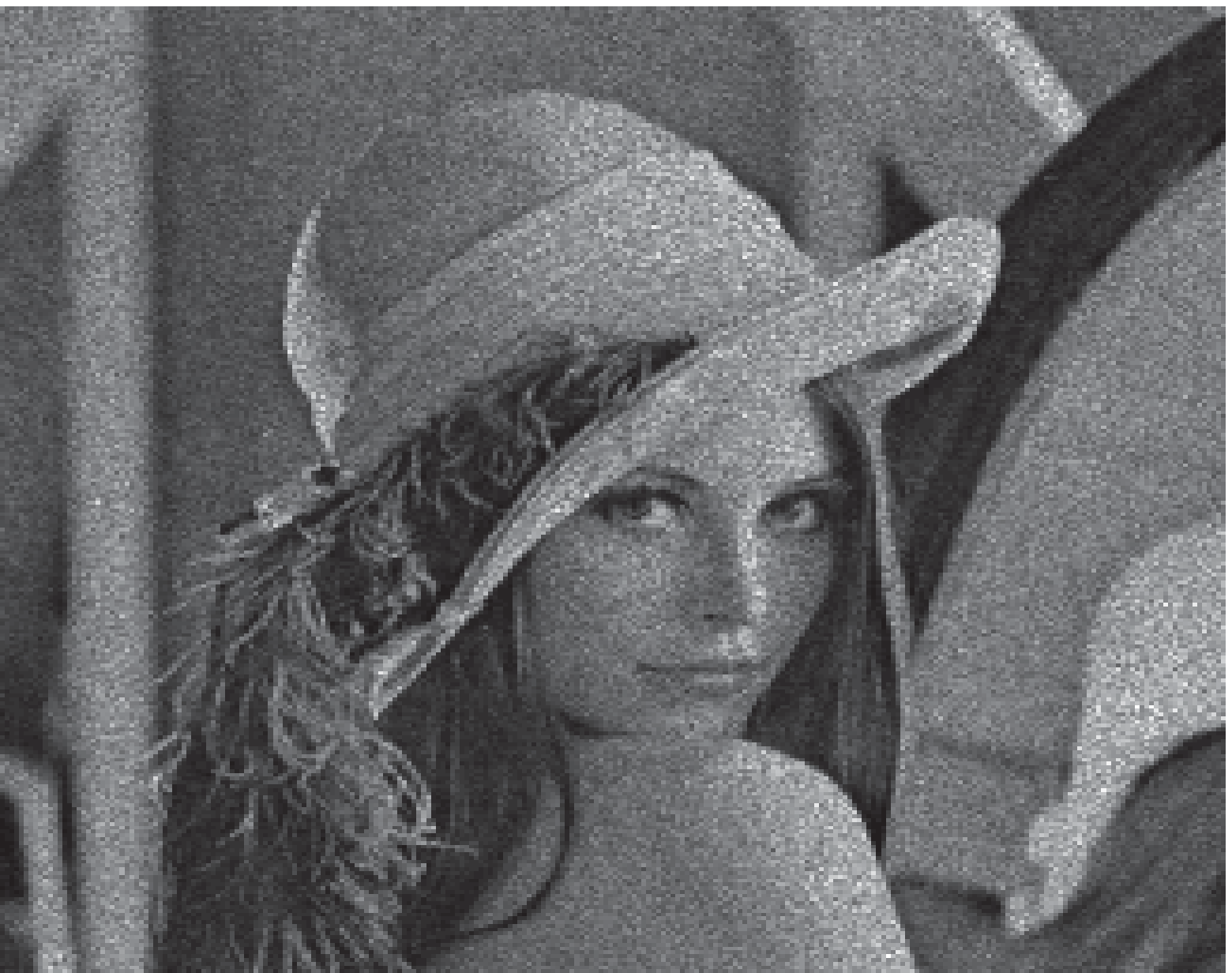}\hspace{0.1cm}\includegraphics[width=4cm]{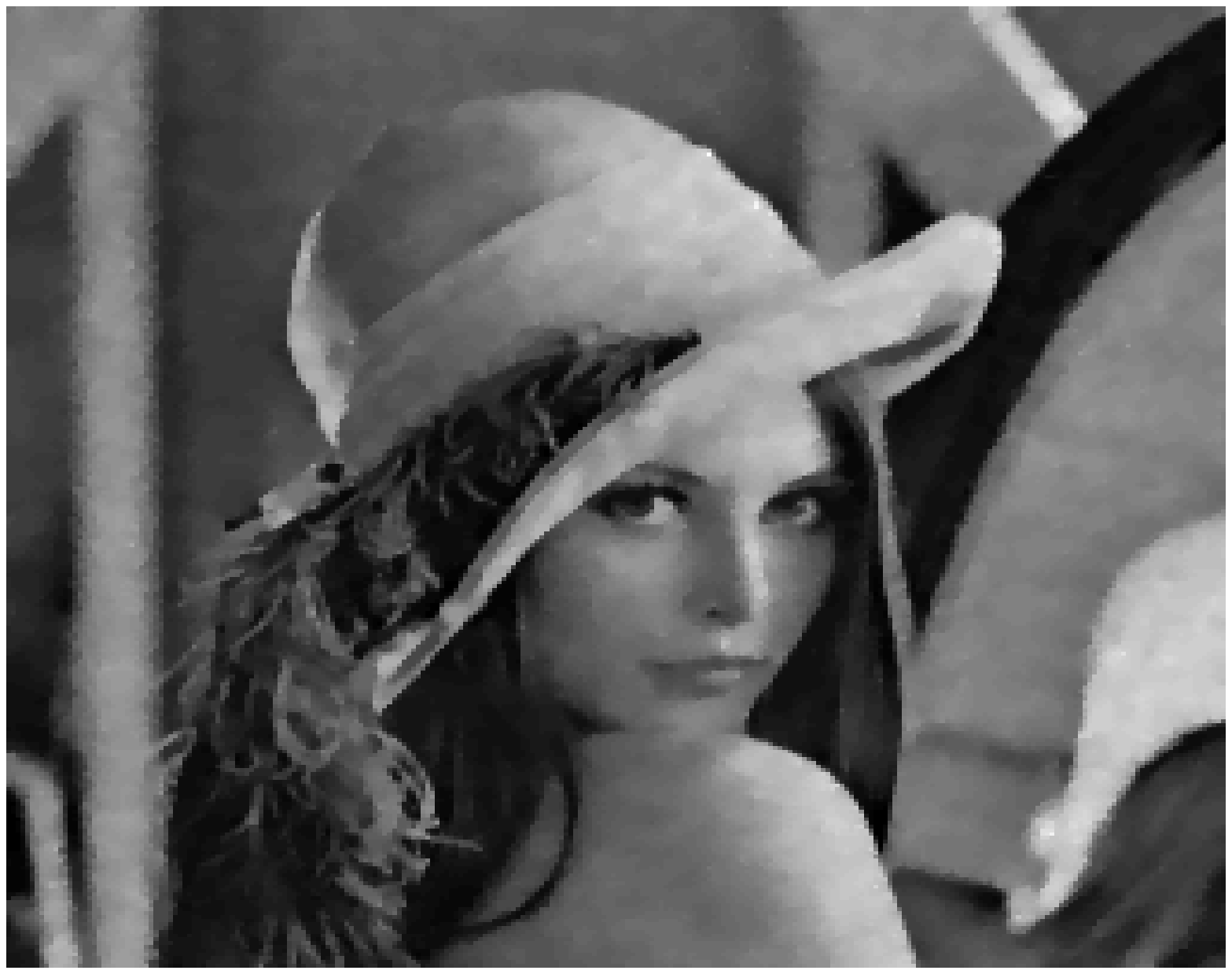}

\vspace{0.1cm}
\includegraphics[width=4cm]{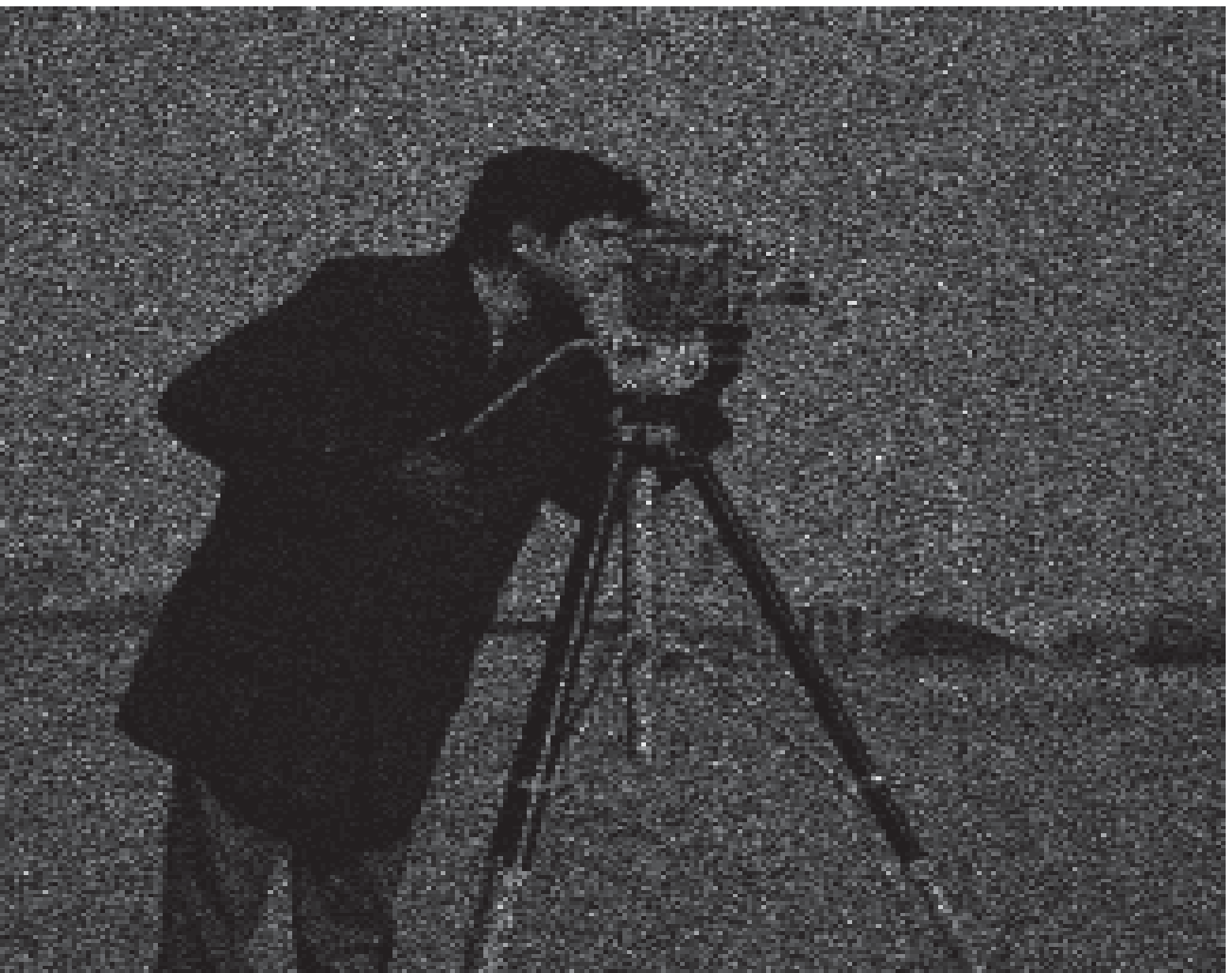}\hspace{0.1cm}\includegraphics[width=4cm]{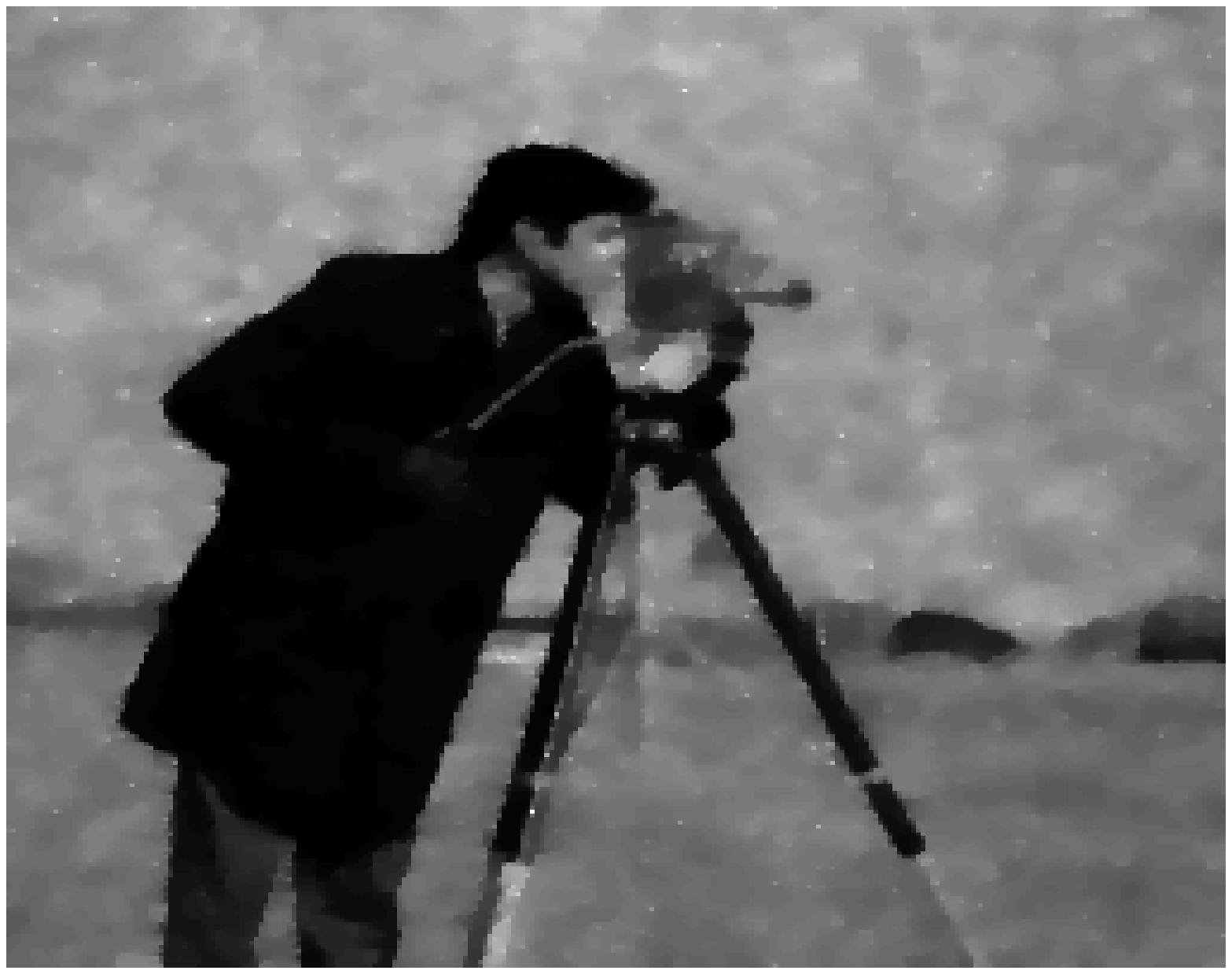}

\vspace{0.1cm}
\includegraphics[width=4cm]{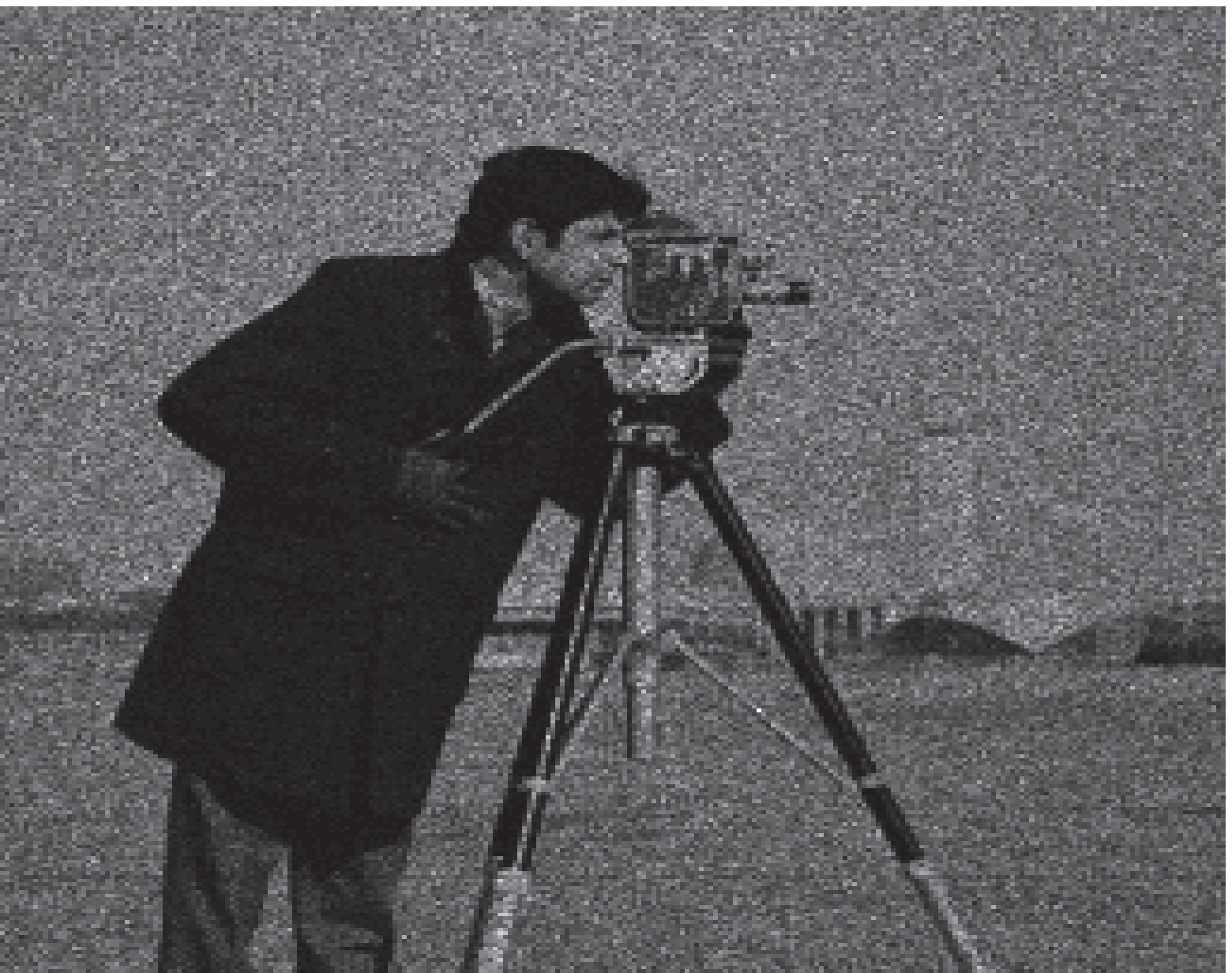}\hspace{0.1cm}\includegraphics[width=4cm]{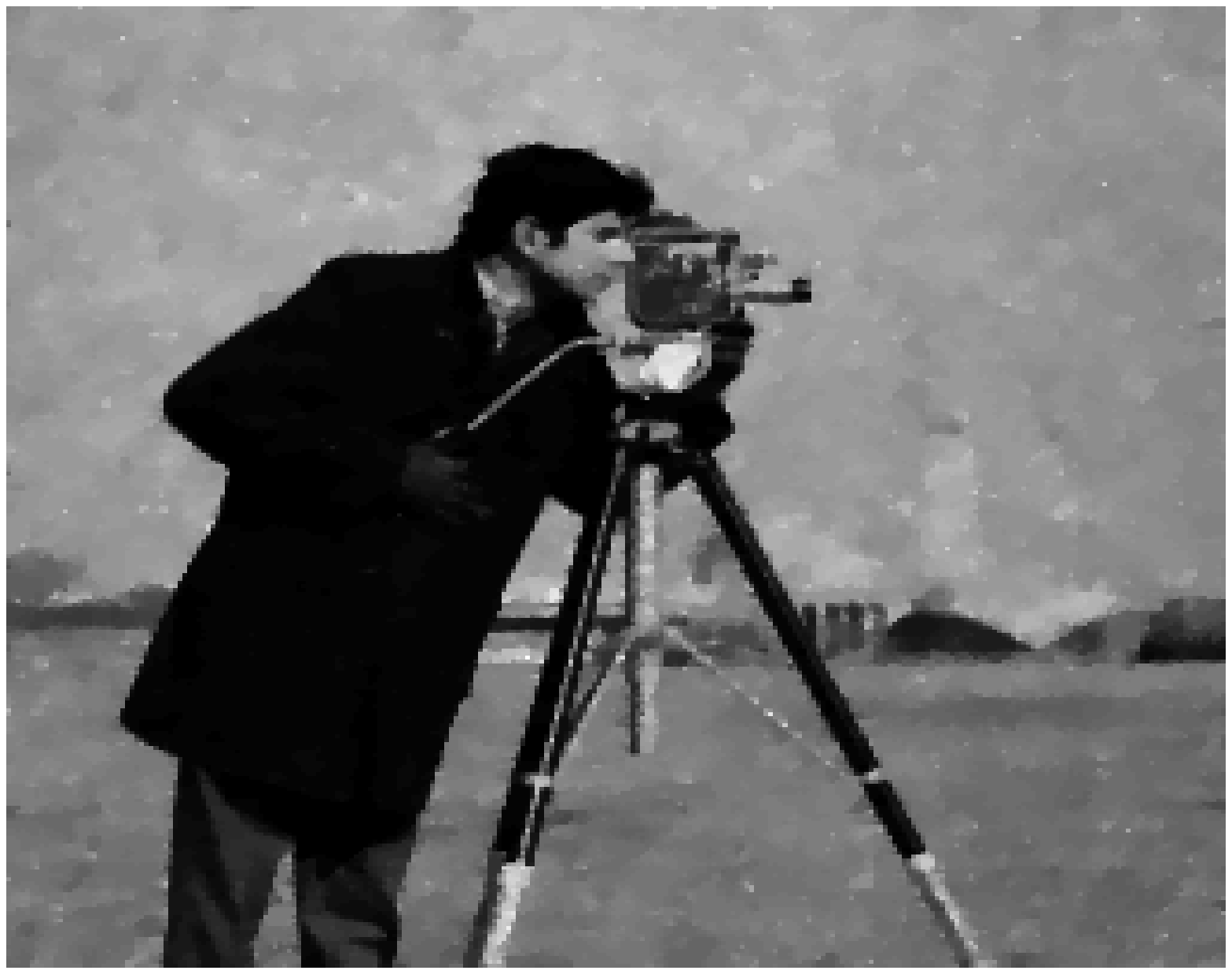}
\caption{Left column: observed noisy images. Right column: image estimates.
First and second rows: Lena, $M=5$ and $M=33$. Third and fourth rows:
Cameraman, $M=3$ and $M=13$.}
\label{fig:images}
\end{figure}

\begin{figure}
\includegraphics[width=4.15cm]{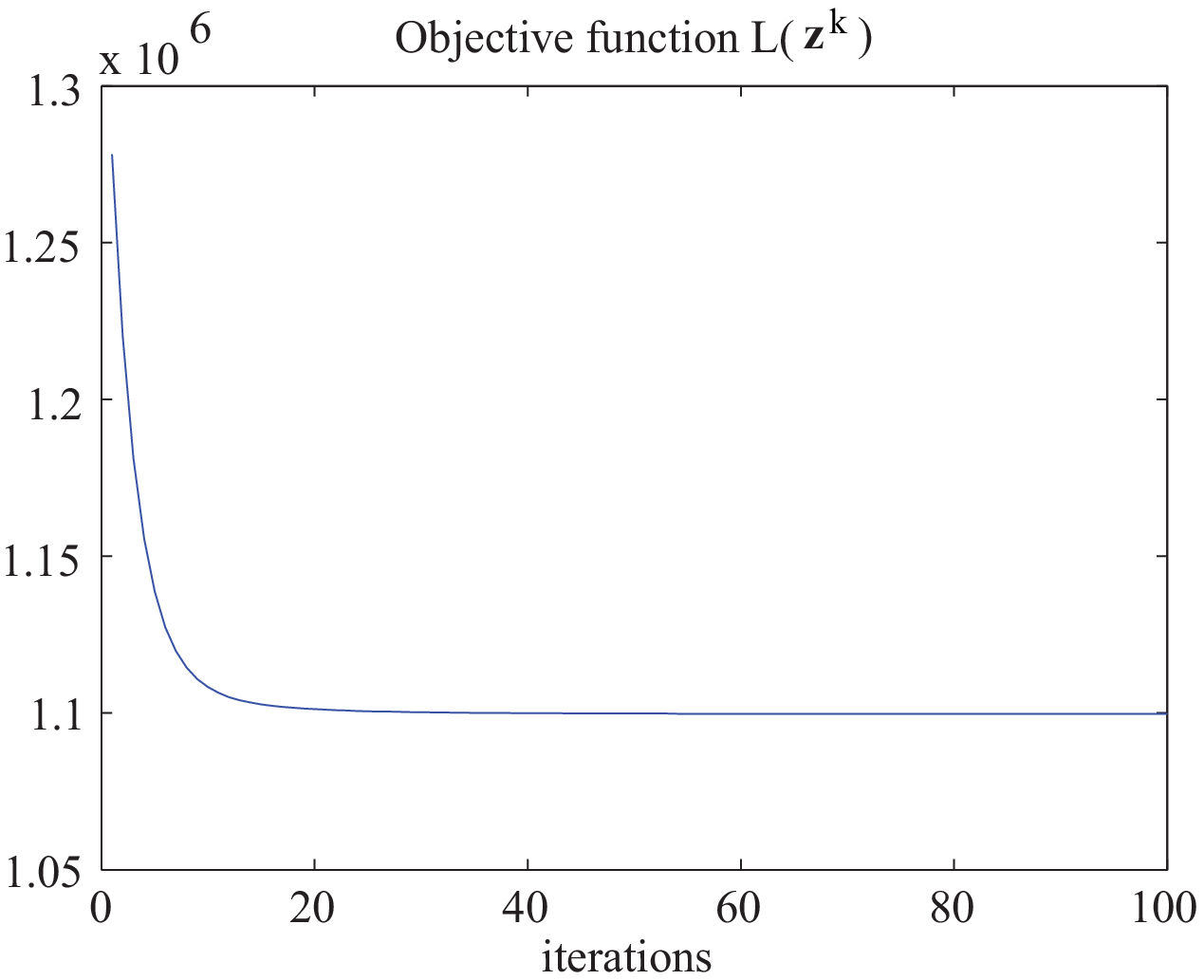}\hspace{.3cm}\includegraphics[width=4.35cm]{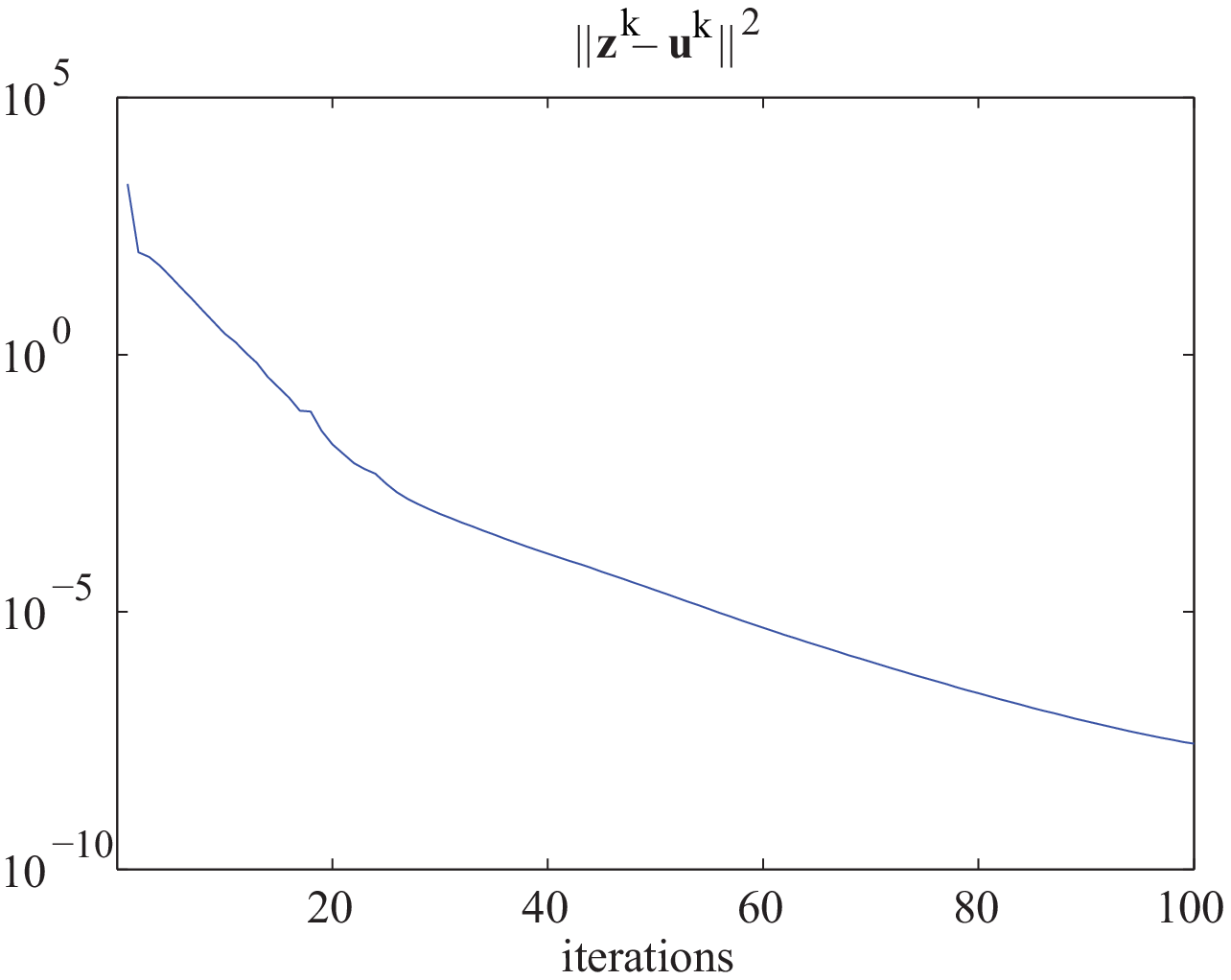}
\caption{Evolution of the objective function $L({\bf z}^k)$ and of the constraint function
$\|{\bf z}^k - {\bf u}^k\|_2^2$, along the iterations of the algorithm, for the
experiment with the Cameraman image and $M=3$. }
\label{fig:plots}
\end{figure}

Figure~\ref{fig:images} shows the noisy and restored images, for the same
experiments reported in Table~\ref{tab:results}. Finally, Figure~\ref{fig:plots}
plots the evolution of the objective function $L({\bf z}^k)$ and of
the constraint function $\|{\bf z}^k - {\bf u}^k\|_2^2$ along the iterations, for
the example with the Cameraman image and $M=3$. Observe the extremely low value of $\|{\bf z}^k - {\bf u}^k\|_2^2$ at the final iterations, showing that, for all practical purposes,
the constraint (\ref{eq:z_u_c}) is satisfied.

\section{Concluding Remarks}
We have proposed an approach to total variation denoising of
images contaminated by multiplicative noise, by exploiting
a split Bregman technique. The proposed algorithm is very
simple and, in the experiments herein reported, exhibited
state of the art performance and speed. We are currently
working on extending our methods to problems involving linear
observation operators ({\it e.g.}, blur) and other related
noise models, such as Poisson.

\end{document}